# Shifted Hankel determinants of Catalan numbers and related results II: Backward shifts

Johann Cigler

johann.cigler@univie.ac.at

**Abstract.**

By prepending zeros to a given sequence Hankel determinants of backward shifts of this sequence become meaningful. We obtain some results for the sequences of Catalan numbers and of some numbers and polynomials which are related to Catalan numbers and propose conjectures for sequences of convolution powers of Catalan numbers.

## 1. Introduction.

In the following all sequences $(a_n)_{n \geq 0}$ will be extended to doubly infinite sequences $(a_n)_{n \in \mathbb{Z}}$ by defining $a_n = 0$ for $n < 0$.

While the Hankel determinants

$$d_m(n) = \det\left(a_{m+i+j}\right)_{i,j=0}^{n-1} \tag{1}$$

of Catalan numbers and related numbers $a_n$ for $m \in \mathbb{N}$ are well known there are apparently no papers dealing with the case $m < 0$.

As is well known, the sequence $(C_n)_{n \geq 0} = (1,1,2,5,14,42,132,429,1430,\cdots)$ of the Catalan numbers $C_n = \frac{1}{n+1}\binom{2n}{n}$ satisfies $d_0(n) = d_1(n) = 1$ and is thereby even uniquely determined.

This fact shows that there is a very strong relationship between Hankel determinants and Catalan numbers. This is reinforced by the fact that for all $m \in \mathbb{N}$ a closed formula for $d_m(n)$ exists.

Let me first recall some well-known results. We shall denote the Hankel determinants of the Catalan numbers by $D_m(n) = \det\left(C_{i+j+m}\right)_{i,j=0}^{n-1}$ with a capital $D$. For $m \in \mathbb{N}$

$$D_m(n) = p_m(n) \tag{2}$$

where $p_m(n)$ is a polynomial of degree $\binom{m}{2}$. For example $p_0(n) = p_1(n) = 1$, $p_2(n) = n+1$,

$$p_3(n) = \frac{1}{6}(1+n)(2+n)(3+2n), \quad p_4(n) = \frac{1}{180}(1+n)(2+n)^2(3+n)(3+2n)(5+2n).$$

For $m \in \mathbb{N}$ $p_m(n)$ has the closed formula

$$p_m(n) = \prod_{1 \leq i \leq j \leq m-1} \frac{2n+i+j}{i+j}. \qquad (3)$$

There are several different proofs of this result (cf. e.g. [5], Theorem 33). An elementary proof can also be found in part I ([3]) of the present article.

It turns out that for $m > 0$

$$D_{-m}(n) = p_{m+1}(-n). \qquad (4)$$

Observing that $p_{m+1}(-n) = \prod_{1 \leq i \leq j \leq m} \frac{-2n+i+j}{i+j} = 0$ for $0 < n \leq m$ by choosing $i = j = n$ for some $n \leq m$ and that for $n \geq m+1$

$$p_{m+1}(-n) = \prod_{1 \leq i \leq j \leq m} \frac{-2n+i+j}{i+j} = \prod_{1 \leq i \leq j \leq m} \frac{-2n+(m+1-i)+(m+1-j)}{i+j} = \prod_{1 \leq i \leq j \leq m} \frac{2(m+1-n)-i-j}{i+j}$$

$$= (-1)^{\binom{m}{2}+m} \prod_{1 \leq i \leq j \leq m} \frac{2(n-m-1)+i+j}{i+j} = (-1)^{\binom{m+1}{2}} p_{m+1}(n-m-1)$$

we get

**Theorem 1**

*If we extend $C_n$ to $n \in \mathbb{Z}$ by setting $C_n = 0$ for $n < 0$, then we get for $m \geq 1$*

$$D_{-m}(n) = \det\left(C_{i+j-m}\right)_{i,j=0}^{n-1} = p_{m+1}(-n) \qquad (5)$$

*or equivalently*

$$D_{-m}(n) = (-1)^{\binom{m+1}{2}} p_{m+1}(n-m-1) = (-1)^{\binom{m+1}{2}} D_{m+1}(n-m-1) \qquad (6)$$

*for $n \geq m+1$ and $D_{-m}(n) = 0$ for $1 \leq n \leq m$.*

For example,

$(D_{-1}(n))_{n \geq 0} = (1, 0, -1, -2, -3, -4, -5, -6, -7, -8, -9, -10, -11, \cdots)$ with

$(D_2(n))_{n \geq 0} = (1, 2, 3, 4, 5, 6, 7, 8, 9, 10, 11, 12, 13, \cdots)$,

$(D_{-2}(n))_{n \geq 0} = (1, 0, 0, -1, -5, -14, -30, -55, -91, -140, -204, -285, -385, \cdots)$ with

$(D_3(n))_{n \geq 0} = (1, 5, 14, 30, 55, 91, 140, 204, 285, 385, \cdots)$,

$(D_{-3}(n))_{n \geq 0} = (1, 0, 0, 0, 1, 14, 84, 330, 1001, 2548, 5712, 11628, 21945, \cdots)$ with

$(D_4(n))_{n \geq 0} = (1, 14, 84, 330, 1001, 2548, 5712, 11628, 21945, 38962, \cdots)$.

Motivated by this fact we obtain similar results about sequences of numbers and polynomials which are in some way related to Catalan numbers.

Let

$$C(x) = \sum_{n \geq 0} C_n x^n = \frac{1 - \sqrt{1-4x}}{2x} \qquad (7)$$

denote the generating function of the Catalan numbers.

We will consider the numbers

$$M_b(n) = \sum_{k=0}^{n} \left( \binom{n+k}{k} - \binom{n+k}{k-1} \right) b^{n-k} \qquad (8)$$

with generating function

$$\sum_{n \geq 0} M_b(n) x^n = \frac{C(x)}{1 - bxC(x)}, \qquad (9)$$

the central binomial coefficients

$$B_n = \binom{2n}{n}, \qquad (10)$$

the Narayana polynomials

$$C_n(t) = \sum_{k=0}^{n} \binom{n-1}{k} \binom{n}{k} \frac{1}{k+1} t^k \qquad (11)$$

which for $t=1$ reduce to $C_n$,

and the Narayana polynomials of typ B

$$B_n(t) = \sum_{k=0}^{n} \binom{n}{k}^2 t^k \qquad (12)$$

which for $t=1$ reduce to $B_n$.

Finally we consider the $k$-th convolution powers of the Catalan numbers

$$C_{k,n} = \binom{2n+k}{n} \frac{k}{2n+k} \qquad (13)$$

with generating function

$$C(x)^k = \sum_{n \geq 0} C_{k,n} x^n \qquad (14)$$

for a positive integer $k$.

## 2. Some background material

### 2.1. Catalan numbers and related numbers

Let us note that $C(x)$ satisfies

$$C(x) = 1 + xC(x)^2 \tag{15}$$

or equivalently

$$C^2(x) = \sum_{n \geq 0} C_{n+1} x^n. \tag{16}$$

The formula (13) for the coefficients of $C(x)^k$ is well known (cf.[4]). For the sake of completeness let us give another simple proof. Let

$$f(k,x) = \sum_{n \geq 0} \left( \binom{2n+k-1}{n} - \binom{2n+k-1}{n-1} \right) x^n \text{ and } g(k,x) = C(x)^k.$$

They satisfy

$$g(k+1,x) = C(x)^{k+1} = C(x)C(x)^k = \left(1 + xC(x)^2\right)C(x)^k = C(x)^k + xC(x)^{k+2} = g(k,x) + xg(k+2,x)$$

and $f(k+1,x) = f(k,x) + xf(k+2,x)$ because

$$\binom{2n+k}{n} - \binom{2n+k}{n-1} = \binom{2n+k-1}{n} - \binom{2n+k-1}{n-1} + \binom{2n-2+k+1}{n-1} - \binom{2n-2+k+1}{n-1-1}.$$

Since $g(1,x) = \sum_{n \geq 0} C_n x^n = \sum_{n \geq 0} \left( \binom{2n}{n} - \binom{2n}{n-1} \right) x^n = f(1,x)$ and

$$g(2,x) = \sum_{n \geq 0} C_{n+1} x^n = \sum_{n \geq 0} \left( \binom{2n+1}{n} - \binom{2n+1}{n-1} \right) x^n = f(2,x)$$

we get by induction $g(k,x) = f(k,x)$ for all $k \in \mathbb{N}$.

For small $k$ we get (cf. OEIS [6], A000108, A000245, A002057)

$$(C_{1,n})_{n \geq 0} = (1,1,2,5,14,42,132,429,1430,4862,\cdots),$$

$$(C_{2,n})_{n \geq 0} = (1,2,5,14,42,132,429,1430,4862,16796,\cdots),$$

$$(C_{3,n})_{n \geq 0} = (1,3,9,28,90,297,1001,3432,11934,41990,\cdots),$$

$$(C_{4,n})_{n \geq 0} = (1,4,14,48,165,572,2002,7072,25194,90440,\cdots).$$

The generating function of the numbers (8) is (9) because

$$\frac{C(x)}{1-bxC(x)} = \sum_{j\geq 0} b^j x^j C(x)^{j+1} = \sum_{j\geq 0} b^j x^j \sum_{k\geq 0}\left(\binom{2k+j}{k}-\binom{2k+j}{k-1}\right)x^k$$

$$= \sum_{n\geq 0} x^n \sum_{j+k=n} b^j \left(\binom{k+n}{k}-\binom{k+n}{k-1}\right).$$

Let us mention some special cases.

1) For $b=0$ we get the Catalan numbers $M_0(n) = C_n$.

2) For $b=1$ we get the shifted Catalan numbers $M_1(n) = C_{n+1}$, because (15) implies

$$\frac{1}{C(x)} = 1 - xC(x) \tag{17}$$

and thus $\dfrac{C(x)}{1-xC(x)} = C(x)^2 = \sum_{n\geq 0} C_{n+1} x^n.$

3) For $b=2$ we get $M_2(n) = \binom{2n+1}{n}$ (cf. OEIS [6], A001700 ),

because $\dfrac{C(x)}{1-2xC(x)} = \dfrac{C(x)}{\sqrt{1-4x}} = \dfrac{1}{2x}\left(\dfrac{1}{\sqrt{1-4x}}-1\right) = \dfrac{1}{2}\sum_{n\geq 0}\binom{2n+2}{n+1}x^n = \sum_{n\geq 0}\binom{2n+1}{n}x^n.$

## 2.2. Hankel determinants

For a sequence $(a_n)_{n\in\mathbb{Z}}$ of real or complex numbers with $a_n = 0$ for $n < 0$ and $a(0) \neq 0$ all Hankel determinants

$$d_m(n) = \det\left(a_{m+i+j}\right)_{i,j=0}^{n-1} \tag{18}$$

for $m \in \mathbb{Z}$ are meaningful for all $n \in \mathbb{N}$. Setting $d_m(0) = 1$ they satisfy

$$\det\begin{pmatrix} d_m(n-1) & d_{m+1}(n-1) \\ d_{m+1}(n-1) & d_{m+2}(n-1) \end{pmatrix} = d_m(n) d_{m+2}(n-2) \tag{19}$$

for $n \geq 2$.

This follows from Dodgson's condensation formula (cf. [5], Proposition 10)

$$\det A \cdot \det A_{1,n}^{1,n} = \det A_1^1 \cdot \det A_n^n - \det A_1^n \cdot \det A_n^1 \tag{20}$$

where $A_{i_1,i_2,\cdots,i_k}^{j_1,j_2,\cdots,j_k}$ denotes the submatrix of a matrix $A$ in which rows $i_1, i_2, \cdots, i_k$ and columns $j_1, j_2, \cdots, j_k$ are omitted.

Since $p_m(n) = D_m(n)$ for the Hankel determinants $D_m(n)$ of the Catalan numbers and $p_m(n) \neq 0$ for $m, n \in \mathbb{N}$ we get

**Lemma 2.**

*For each $n \in \mathbb{N}$ the numbers $p_m(n) = \prod_{1 \leq i \leq j \leq m-1} \dfrac{2n+i+j}{i+j}$, $m \in \mathbb{N}$, are uniquely determined by the recursion*

$$\det\begin{pmatrix} p_m(n-1) & p_{m+1}(n-1) \\ p_{m+1}(n-1) & p_{m+2}(n-1) \end{pmatrix} = p_m(n) p_{m+2}(n-2) \tag{21}$$

*with the initial values $p_m(0) = 1$ and $p_m(1) = D_m(1) = C_m$ for all $m \in \mathbb{N}$.*

We are especially interested in Hankel matrices of the form $\mathbf{A}_{-m}(n) = \left(a_{i+j-m}\right)_{i,j=0}^{n-1}$ for $m \in \mathbb{N}$. Often it will be convenient to write instead

$$\mathbf{V}_k(n) = \mathbf{A}_{k-n}(n). \tag{22}$$

Then for $k \leq n$ the first row of $\mathbf{V}_k(n)$ contains precisely $k$ entries $a_j$ with $j \geq 0$.

For example,

$$\mathbf{A}_{-3}(4) = \mathbf{V}_1(4) = \begin{pmatrix} 0 & 0 & 0 & a_0 \\ 0 & 0 & a_0 & a_1 \\ 0 & a_0 & a_1 & a_2 \\ a_0 & a_1 & a_2 & a_3 \end{pmatrix}, \quad \mathbf{A}_{-2}(4) = \mathbf{V}_2(4) = \begin{pmatrix} 0 & 0 & a_0 & a_1 \\ 0 & a_0 & a_1 & a_2 \\ a_0 & a_1 & a_2 & a_3 \\ a_1 & a_2 & a_3 & a_4 \end{pmatrix}.$$

Setting

$$v_k(n) = \det \mathbf{V}_k(n) \tag{23}$$

it is obvious that

$$v_1(n) = (-1)^{\binom{n}{2}} a_0^n. \tag{24}$$

For $n \geq 2$ $\mathbf{V}_2(n)$ is obtained from $\mathbf{V}_1(n+1)$ by deleting the first row and column. By Cramer's rule

$\left(\mathbf{V}_1(n+1)\right)^{-1} = \dfrac{1}{\det\left(\mathbf{V}_1(n+1)\right)} \left(\alpha_{j,i}\right)_{i,j=0}^n$ with $\alpha_{j,i} = (-1)^{i+j} \det A_{j,i}$, where $A_{i,j}$ is the matrix obtained by crossing out row $i$ and column $j$ in $\mathbf{V}_1(n+1)$. Thus $A_{0,0} = \mathbf{V}_2(n)$.

Therefore $(-1)^{\binom{n+1}{2}} \dfrac{v_2(n)}{a_0^{n+1}}$ is the entry in position $(0,0)$ of the inverse matrix of $\mathbf{V}_1(n+1)$.

To obtain this inverse we use

**Lemma 3**

Let $a(n) = b(n) = c(n) = 0$ for $n < 0$. Then

$$\left(a(i+j-n)\right)_{i,j=0}^{n} \left(b(n-j-k)\right)_{j,k=0}^{n} = \left(c(i-k)\right)_{i,k=0}^{n} \text{ if and only if}$$

$$\sum_{n \geq 0} a(n)x^n \sum_{n \geq 0} b(n)x^n = \sum_{n \geq 0} c(n)x^n.$$

For example

$$\begin{pmatrix} 0 & 0 & 0 & a(0) \\ 0 & 0 & a(0) & a(1) \\ 0 & a(0) & a(1) & a(2) \\ a(0) & a(1) & a(2) & a(3) \end{pmatrix} \begin{pmatrix} b(3) & b(2) & b(1) & b(0) \\ b(2) & b(1) & b(0) & 0 \\ b(1) & b(0) & 0 & 0 \\ b(0) & 0 & 0 & 0 \end{pmatrix} = \begin{pmatrix} c(0) & 0 & 0 & 0 \\ c(1) & c(0) & 0 & 0 \\ c(2) & c(1) & c(0) & 0 \\ c(3) & c(2) & c(1) & c(0) \end{pmatrix}.$$

The proof is obvious because

$$\sum a(i+j-n)b(n-j-k) = \sum_{j \geq n-i} a(j+i-n)b(n-j-k) = \sum_{j \geq 0} a(j)b(i-k-j) = c(i-k).$$

**Corollary 4**

If $\dfrac{1}{\sum_{n \geq 0} a(n)x^n} = \sum_{n \geq 0} b(n)x^n$ then

$$v_2(n) = \det\left(a(2-n+i+j)\right)_{i,j=0}^{n-1} = (-1)^{\binom{n+1}{2}} a(0)^{n+1} b(n). \tag{25}$$

### 3. Computation of the Hankel determinants

Dodgson's condensation formula (20) gives

$$v_k(n+k)v_k(n+k-2) = v_{k-1}(n+k-1)v_{k+1}(n+k-1) - v_k(n+k-1)^2 \tag{26}$$

for $n \geq k$.

**Lemma 5**

Let $(a_n)_{n \in \mathbb{Z}}$ with $a_n = 0$ for $n < 0$ and $a_0 = 1$. If $v_2(n) = (-1)^{\binom{n-1}{2}} aC_{n-1}$ then

$$v_k(n+k) = (-1)^{\binom{n+1}{2}} a^{k-1} p_{n+1}(k-1). \tag{27}$$

**Proof**

Let $v_k(n) = (-1)^{\binom{n-k+1}{2}} a^{k-1} r(n-k+1, k-1)$ then (26) becomes

$$\det \begin{pmatrix} a^{k-1} r(n-1, k-1) & a^{k-1} r(n, k-1) \\ a^{k-1} r(n, k-1) & a^{k-1} r(n+1, k-1) \end{pmatrix} = a^{2k-2} r(n+1, k-2) r(n-1, k)$$

and by cancelling powers of $a$

$$\det\begin{pmatrix} r(n-1,k-1) & r(n,k-1) \\ r(n,k-1) & r(n+1,k-1) \end{pmatrix} = r(n+1,k-2)r(n-1,k)$$

with $r(n,0)=1$ and $r(n-1,1)=C_{n-1}$. Note that $(-1)^{\binom{n}{2}+\binom{n+2}{2}}=-1$.

By Lemma 2 we get $r(n,k)=p_n(k)$ and therefore (27).

### 3.1. The numbers $M_b(n)$.

**Theorem 6**

Let

$$D_m(b;n) = \det\left(M_b(m+i+j)\right)_{i,j=0}^{n-1} \tag{28}$$

denote the Hankel determinants of the numbers $M_b(n)$.

If we extend $M_b(n)$ to $n \in \mathbb{Z}$ by setting $M_b(n)=0$ for $n<0$, then we get for $m \geq 1$

$$D_{-m}(b;n) = \det\left(M_b(-m+i+j)\right)_{i,j=0}^{n-1} = p_{m+1}(-n) \tag{29}$$

or equivalently

$$D_{-m}(b;n) = (-1)^{\binom{m+1}{2}} p_{m+1}(n-m-1) \tag{30}$$

for $n \geq m+1$ and $D_{-m}(b;n)=0$ for $1 \leq n \leq m$.

**Proof**

Using (17) we get

$$\frac{1}{\sum_{n\geq 0} M_b(n)x^n} = \frac{1-bxC(x)}{C(x)} = \frac{1}{C(x)} - bx = 1 - xC(x) - bx = 1-(1+b)x - \sum_{n\geq 2} C_{n-1}x^n.$$

Corollary 4 gives $v_2(b;n) = (-1)^{\binom{n-1}{2}} C_{n-1}$ and Lemma 5 $v_k(n+k) = (-1)^{\binom{n+1}{2}} p_{n+1}(k-1)$.

Therefore, $D_{-m}(b;n) = V_{n-m}(n) = (-1)^{\binom{m+1}{2}} p_{m+1}(n-m-1)$.

The special case $b=0$ gives Theorem 1.

For example for $b=0$ and $b=1$ we get

$$\det\begin{pmatrix} 0 & 1 & 1 & 2 \\ 1 & 1 & 2 & 5 \\ 1 & 2 & 5 & 14 \\ 2 & 5 & 14 & 42 \end{pmatrix} = \det\begin{pmatrix} 0 & 1 & 2 & 5 \\ 1 & 2 & 5 & 14 \\ 2 & 5 & 14 & 42 \\ 5 & 14 & 42 & 132 \end{pmatrix} = -3 = (-1)^{\binom{2}{2}} p_2(2) = p_2(-4),$$

$$\det\begin{pmatrix} 0 & 0 & 1 & 1 \\ 0 & 1 & 1 & 2 \\ 1 & 1 & 2 & 5 \\ 1 & 2 & 5 & 14 \end{pmatrix} = \det\begin{pmatrix} 0 & 0 & 1 & 2 \\ 0 & 1 & 2 & 5 \\ 1 & 2 & 5 & 14 \\ 2 & 5 & 14 & 42 \end{pmatrix} = -5 = (-1)^{\binom{3}{2}} p_3(1) = p_3(-4).$$

### 3.2. Central binomial coefficients

**Theorem 7**

Let $B_n = \binom{2n}{n}$ and $B_n = 0$ for $n < 0$ and let $D_m(B;n) = \det\left(B_{i+j+m}\right)_{i,j=0}^{n-1}$. Then we get for $m \geq 1$

$$D_{-m}(B;n) = (-1)^{\binom{m+1}{2}} 2^{n-m-1} p_{m+1}(n-m-1) \tag{31}$$

for $n \geq m+1$ and $D_{-m}(B;n) = 0$ for $1 \leq n \leq m$.

**Proof**

Here we have $\sum_{n \geq 0} \binom{2n}{n} x^n = \frac{1}{\sqrt{1-4x}}$ and therefore $\frac{1}{\sum_{n \geq 0} \binom{2n}{n} x^n} = \sqrt{1-4x} = 1 - 2xC(x).$

### 3.3. Narayana polynomials

Let

$$\left(C_n(t)\right)_{n \geq 0} = \left(1, 1, 1+t, 1+3t+t^2, 1+6t+6t^2+t^3, 1+10t+20t^2+10t^3+t^4, \cdots\right)$$

be the Narayana polynomials

$$C_n(t) = \sum_{k=0}^{n} \binom{n-1}{k}\binom{n}{k}\frac{1}{k+1} t^k. \tag{32}$$

Their Hankel determinants $D_m(C(t);n) = \det\left(C_{i+j+m}(t)\right)_{i,j=0}^{n-1}$ show for negative $m$ a similar behavior as their counterparts for the Catalan numbers.

**Theorem 8**

*If we extend $C_n(t)$ to $n \in \mathbb{Z}$ by setting $C_n(t) = 0$ for $n < 0$, then we get for $m \geq 1$*

$$D_{-m}(C(t);n) = (-1)^{\binom{m+1}{2}} t^{n-m-1} D_{m+1}(C(t); n-m-1) \tag{33}$$

*for $n \geq m+1$ and $D_{-m}(C(t);n) = 0$ for $1 \leq n \leq m$.*

**Proof**

The generating function $C(t,x) = \sum_{n \geq 0} C_n(t) x^n$ of the Narayana polynomials is (cf. [7])

$$C(t,x) = \frac{1 + x(t-1) - \sqrt{1 - 2x(t+1) + x^2(t-1)^2}}{2tx}. \tag{34}$$

It satisfies (cf. [7])

$$C(t,x) = 1 + (1-t)xC(t,x) + txC(t,x)^2. \tag{35}$$

Therefore

$$\frac{1}{C(t,x)} = 1 + (t-1)x - txC(t,x). \tag{36}$$

Here we need a modification of Lemma 2:

*For each $n \in \mathbb{N}$ the polynomials $p_m(t,n) = D_m(C(t);n)$ $m \in \mathbb{N}$, are uniquely determined by the recursion*

$$\det \begin{pmatrix} p_m(t,n-1) & p_{m+1}(t,n-1) \\ p_{m+1}(t,n-1) & p_{m+2}(t,n-1) \end{pmatrix} = p_m(t,n) p_{m+2}(t,n-2) \tag{37}$$

*with the initial values $p_m(t,0) = 1$ and $p_m(t,1) = D_m(C(t),1) = C_m(t)$ for all $m \in \mathbb{N}$.*

Note that $p_m(t,n) = D_m(C(t);n)$ as a determinant with polynomial entries is also a polynomial in $t$. No $p_m(t,n)$ vanishes because $p_m(1,n) = p_m(n) \neq 0$.

By Corollary 4 and Lemma 5 we get (33).

The first polynomials are $p_0(t,n) = p_1(t,n) = t^{\binom{n}{2}}$, $p_2(t,n) = t^{\binom{n}{2}}(1 + t + \cdots + t^n) = t^{\binom{n}{2}} \frac{1-t^{n+1}}{1-t}$.

Apparently, there is no closed formula known for the general case.

For example

$$(D_{-1}(C(t);n))_{n \geq 0} = \left(1, 0, -1, -t(1+t), -t^3(1+t+t^2), -t^6(1+t)(1+t^2), -t^{10}(1+t+t^2+t^3+t^4), \cdots\right)$$

with

$$(D_2(C(t);n))_{n\geq 0} = \left(1, 1+t, t(1+t+t^2), t^3(1+t)(1+t^2), t^6(1+t+t^2+t^3+t^4), t^{10}(1+t)(1-t+t^2)(1+t+t^2), \cdots\right),$$

$$(D_{-2}(C(t);n))_{n\geq 0} = \left(1, 0, 0, -1, -t(1+3t+t^2), -t^3(1+3t+6t^2+3t^3+t^4), -t^6(1+3t+6t^2+10t^3+6t^4+3t^5+t^6), \cdots\right)$$

with

$$(D_3(C(t);n))_{n\geq 0} = \left(1, 1+3t+t^2, t(1+3t+6t^2+3t^3+t^4), t^3(1+3t+6t^2+10t^3+6t^4+3t^5+t^6), \cdots\right),$$

$$(D_{-3}(C(t);n))_{n\geq 0} = \left(1, 0, 0, 0, 1, t+6t^2+6t^3+t^4, t^3+6t^4+21t^5+28t^6+21t^7+6t^8+t^9, \cdots\right)$$

with

$$(D_4(C(t);n))_{n\geq 0} = \left(1, 1+6t+6t^2+t^3, t+6t^2+21t^3+28t^4+21t^5+6t^6+t^7, \cdots\right).$$

### 3.4. Narayana polynomials of type B

Let $B_n(t) = \sum_{k=0}^{n} \binom{n}{k}^2 t^k$ be the Narayana polynomials of type B which for $t=1$ reduce to the central binomial coefficients and denote by $D_m(B(t), n)$ their Hankel determinants.

**Theorem 9**

Let $B_n(t) = \sum_{k=0}^{n} \binom{n}{k}^2 t^k$ and $B_n(t) = 0$ for $n < 0$. Then we get for $m \geq 1$

$$D_{-m}(B(t);n) = (-1)^{\binom{m+1}{2}} (2t)^{n-m-1} p_{m+1}(t, n-m-1) \tag{38}$$

for $n \geq m+1$ and $D_{-m}(B(t);n) = 0$ for $1 \leq n \leq m$.

**Proof**

The generating function of $B_n(t)$ is (cf. [3])

$$\sum_{n\geq 0} B_n(t) x^n = \frac{1}{\sqrt{(1-(1+t)x)^2 - 4tx^2}} = \frac{1}{1+(t-1)x - 2txC(t,x)}$$

The result follows as above for $C(t)$. For example,

$$(D_{-1}(B(t);n))_{n\geq 0} = \left(1, 0, -1, -2t(1+t), -4t^3(1+t+t^2), -8t^6(1+t)(1+t^2), -16t^{10}(1+t+t^2+t^3+t^4), \cdots\right),$$

$$(D_{-2}(B(t);n))_{n\geq 0} = \left(1, 0, 0, -1, -2t(1+3t+t^2), -4t^3(1+3t+6t^2+3t^3+t^4), \cdots\right),$$

$$(D_{-3}(B(t);n))_{n\geq 0} = \left(1, 0, 0, 0, 1, 2t(1+t)(1+5t+t^2), 4t^3(1+6t+21t^2+28t^3+21t^4+6t^5+t^6), \cdots\right).$$

### 3.5. Convolution powers of Catalan numbers

Finally consider the $k$-th convolution powers of the Catalan numbers $C_{k,n} = \binom{2n+k}{n}\frac{k}{2n+k}$ with generating function $C(x)^k = \sum_{n \geq 0} C_{k,n} x^n$ for a positive integer $k$ and $n \in \mathbb{N}$.

We extend $C_{k,n}$ to all $n \in \mathbb{Z}$ by setting $C_{k,n} = 0$ for $n < 0$.

Let

$$D_{k,m}(n) = \det\left(C_{k,i+j+m}\right)_{i,j=0}^{n-1} \tag{39}$$

be their Hankel determinants.

For $k > 2$ little is known. Computations suggest that these determinants show remarkable modular patterns (cf. [2], [8]). For example, it seems that

$D_{3,0}(3n) = D_{3,0}(3n+1) = (-1)^n, \quad D_{3,0}(3n+2) = 0.$

$D_{4,0}(2n) = D_{4,0}(2n+1) = (-1)^n (n+1),$

$D_{5,0}(5n) = D_{5,0}(5n+1) = 1, \quad D_{5,0}(5n+2) = -5(n+1), \quad D_{5,0}(5n+3) = 0, \quad D_{5,0}(5n+4) = 5(n+1),$

$D_{6,0}(3n) = D_{6,0}(3n+1) = (-1)^n (n+1)^2, \quad D_{6,0}(3n+2) = (-1)^{n+1}\frac{3}{2}(1+n)(2+n)(3+2n),$

$D_{7,0}(7n) = D_{7,0}(7n+1) = (-1)^n, \quad D_{7,0}(7n+2) = (-1)^n \frac{7}{6}(1+n)(-12+49n+98n^2),$

$D_{7,0}(7n+3) = (-1)^{n+1} 7^2 (n+1)^2, \quad D_{7,0}(7n+4) = 0, \quad D_{7,0}(7n+5) = (-1)^n 7^2 (n+1)^2,$

$D_{7,0}(7n+6) = (-1)^n \frac{7}{6}(1+n)(282+343n+98n^2).$

Interestingly, a very general generalization of Theorem 1 seems to hold for these sequences:

**Conjecture 10**

*For a positive integer $k$ let $C_{k,n} = \binom{2n+k}{n}\frac{k}{2n+k}$ and $C_{k,n} = 0$ for $n < 0$.*

*Then we get for $m \in \mathbb{N}$*

$$D_{2k,1-k-m}(n) = (-1)^{\binom{m+k}{2}} D_{2k,1-k+m}(n-m-k) \tag{40}$$

*for $n \geq m+k$ and $D_{2k,1-k-m}(n) = 0$ for $0 < n < m+k$*

*and*

$$D_{2k-1,2-k-m}(n) = (-1)^{\binom{m+k-1}{2}} D_{2k-1,-k+m+1}(n-m-k+1) \quad (41)$$

for $n \geq m+k-1$ and $D_{2k-1,2-k-m}(n) = 0$ for $0 < n < m+k-1$.

For $m = 0$ this implies

**Conjecture 11**

$$D_{2k,1-k}(kn) = (-1)^{\binom{k}{2}n} \quad (42)$$
$$D_{2k,1-k}(n) = 0 \quad \text{else.}$$

and

$$D_{2k-1,2-k}((2k-1)n) = (-1)^{(k-1)n},$$
$$D_{2k-1,2-k}((2k-1)n+k-1) = (-1)^{(k-1)n+\binom{k-1}{2}}, \quad (43)$$
$$D_{2k-1,2-k}(n) = 0 \quad \text{else.}$$

**Proof**

By (40) we get $D_{2k,1-k}(n) = (-1)^{\binom{k}{2}} D_{2k,1-k}(n-k)$ and $D_{2k,1-k}(n) = 0$ for $1 < n < k$. Therefore, $D_{2k,1-k}(0) = 1$, $D_{2k,1-k}(k) = (-1)^{\binom{k}{2}}$ and $D_{2k,1-k}(n) = 0$ for $0 < n < k$, which gives (42).

By (41) we get for $m = 0$ $D_{2k-1,2-k}(n) = (-1)^{\binom{k-1}{2}} D_{2k-1,-k+1}(n-k+1)$ and $D_{2k-1,2-k}(n) = 0$ for $0 < n < k$ and for $m = 1$ $D_{2k-1,1-k}(n) = (-1)^{\binom{k}{2}} D_{2k-1,-k+2}(n-k)$ for $0 < n < k+1$.

Therefore $D_{2k-1,2-k}(0) = 1$, $D_{2k-1,2-k}(k-1) = (-1)^{\binom{k-1}{2}}$, $D_{2k-1,2-k}(2k-1) = (-1)^{\binom{k-1}{2}+\binom{k}{2}} = (-1)^{k-1}$, and else $0$ in $\{0, 2k-1\}$. This implies (43).

Some examples for Conjecture 10:

If we set $d_{2k,m}(n) = D_{2k,m+1-k}(n)$ and $d_{2k-1,m} = D_{2k-1,m+2-k}(n)$ then

$$(d_{3,-2}(n))_{n \geq 0} = (1,0,0,-1,-3,-3,1,6,6,-1,-9,-9,1,12,12,\cdots),$$

$$(d_{3,1}(n))_{n \geq 0} = (1,3,3,-1,-6,-6,1,9,9,-1,-12,-12,1,15,15,\cdots),$$

$$(d_{4,-2}(n))_{n \geq 0} = (1,0,0,0,1,4,-4,-20,9,56,-16,-120,25,220,-36,\cdots)$$

$$(d_{4,2}(n)) = (1,4,-4,-20,9,56,-16,-120,25,220,-36,-364,49,560,-64,\cdots).$$

Finally, I want to mention a Conjecture which could be interpreted as a generalization of the identities $D_1(n) = 1$ and $D_2(n) = n+1$ to which it reduces for $k = 1$.

**Conjecture 12**

$$D_{2k,m+1-k}(kn) = (-1)^{\binom{k}{2}n}(n+1)^m,$$

$$D_{2k-1,m+2-k}((2k-1)n+k-1) = (-1)^{\binom{k-1}{2}+(k-1)n}(2k-1)^m(n+1)^m \quad (44)$$

for $0 \le m \le k$.